\documentclass[11pt]{article}
\usepackage{amsmath,amssymb}
\usepackage{graphicx}
\usepackage{mathrsfs}
\usepackage{bbm}
\usepackage{xcolor}
\usepackage{stmaryrd}
\usepackage{subfigure}
\usepackage{tikz}
\usetikzlibrary{arrows}
\tikzstyle{block}=[draw opacity=0.7,line width=1.4cm]

\newtheorem{remark}{Remark}
\newtheorem{theorem}{Theorem}

\newtheorem{definition}{Definition}
\newtheorem{corollary}{Corollary}

\def \R{{\mathbb{R}}}
\def \N{{\mathbb{N}}}

\title{A constructive
approach to regularity of 
Lagrangian trajectories for 
incompressible Euler flow in a bounded domain}

\author{ Nicolas Besse  \footnotemark[2]\ , Uriel Frisch \footnotemark[3]} 

\begin{document}

\maketitle

\renewcommand{\thefootnote}{\fnsymbol{footnote}}

\footnotetext[2]{ 
 Laboratoire J.-L.~Lagrange,  UCA, OCA,
 CNRS, CS 34229, 06304 Nice Cedex 4 France
({\tt Nicolas.Besse@oca.eu})}
  
\footnotetext[3]{
Laboratoire J.-L.~Lagrange,  UCA, OCA,
 CNRS, CS 34229, 06304 Nice Cedex 4 France
({\tt uriel@oca.eu})}

\renewcommand{\thefootnote}{\arabic{footnote}}

\begin{abstract}

The 3D incompressible Euler equation is an important research topic in
the mathematical study of fluid dynamics.  Not only is the global
regularity for smooth initial data an open issue, but the behaviour
may also depend on the presence or absence of boundaries.

For a good understanding, it is crucial to carry out, besides
mathematical studies, high-accuracy and well-resolved numerical
exploration.  Such studies can be very demanding in computational
resources, but recently it has been shown that very substantial gains
can be achieved first, by using Cauchy's Lagrangian formulation of the
Euler equations and second, by taking advantages of analyticity
results of the Lagrangian trajectories for flows whose initial
vorticity is H\"older-continuous.  The latter has been known for about
twenty years (Serfati, J. Math. Pures Appl.  74: 95-104, 1995), but the
combination of the two, which makes use of recursion relations among
time-Taylor coefficients to obtain constructively the time-Taylor
series of the Lagrangian map, has been achieved only recently (Frisch
and Zheligovsky,   Commun. Math. Phys. 326: 499--505, 2014; Podvigina
{\em et al.},  J. Comput. Phys. 306: 320-342, 2016 and references therein).

Here we extend this methodology to incompressible Euler flow in an
impermeable bounded domain whose boundary may be either analytic
or have a regularity between indefinite differentiability and
analyticity.  Non-constructive regularity results for these cases have
already been obtained by Glass {\em et al.}
(Ann. Scient. \'Ec. Norm. Sup. $4^{e}$ s\'erie 45: 1-51, 2012). Using the
invariance of the boundary under the Lagrangian flow, we
establish novel recursion relations that include contributions from
the boundary. This leads to a constructive proof of time-analyticity
of the Lagrangian trajectories with analytic boundaries, which can
then be used subsequently for the design of a very high-order
Cauchy--Lagrangian method.

\end{abstract}

\begin{flushleft}
\end{flushleft}

\begin{flushleft}
{\bf AMS: } to be completed
\end{flushleft}

\section{Introduction}
\label{sec:intro}
The Cauchy problem for the 3D incompressible Euler equation is known to be well posed in time
when the initial velocity is in H\"older or Sobolev spaces with suitable indexes of regularity.
This was established in the seminal work of  Lichtenstein \cite{Lic27, Lic29} and Gyunter \cite{Gun26, Gun34} 
for the case of the whole space and, of Ebin and Marsden \cite{EM70} for the case of bounded domains.
In other words, limited initial
spatial smoothness is preserved in time, at least for a finite duration. Here is a typical result: Suppose the initial
velocity $v_0\in \mathscr{C}^{\,m,\gamma}$, with $m\in \N-\{0\}$ and $0<\gamma<1$
(see, e.g., \cite{Ada75} for the definition of such H\"older spaces).
It needs to be just a little better than  $\mathscr{C}^{1}$ for the
existence  of  a unique solution $v(t,\cdot) \in \mathscr{C}^{\,m,\gamma}$ for $|t|<T(\|v_0\|_{\mathscr{C}^{\,m,\gamma}})$. 
Moreover,  the solution $v(t,\cdot)$
has exactly the same spatial regularity as the initial velocity $v_0$, i.e. if
$v(t,\cdot) \in \mathscr{C}^{\,n,\eta}$ with $n>m$ and $\eta>\gamma$ for some $t>0$, then
 $v_0\in \mathscr{C}^{\,n,\eta}$. Nevertheless, the issue of loss of regularity after a finite time (blow-up) 
of solutions for the 3D incompressible Euler equation still remains an open problem.

In Eulerian coordinates, lack of regularity in space also translates into a lack of regularity in time.
For instance, in Ref. \cite{BDLIS15} are constructed periodic weak solutions 
of the 3D incompressible Euler equation, which dissipate the
total energy and which belong to the H\"older space of exponent $\gamma < 1/5$ ($m=0$), in both space and time.
By contrast in Lagrangian coordinates, where one is focusing on particle trajectories, 
an initial velocity field with limited smoothness (typically in H\"older classes with suitable indexes of regularity) 
launches characteristic curves whose temporal smoothness (typically in some indefinitely differentiable function classes)
widely exceeds the limited spatial smoothness. These are the kinds of typical statements that we will prove here.

Actually, it was already pointed out in Ref. \cite{EM70} that the Lagrangian structure 
can be nicer than the Eulerian one. For this, it is convenient to use the geometric interpretation 
\cite{Arn66} of the Euler equations in Lagrangian coordinates to show that 
``no loss of derivative occurs''. From a technical point of view, the simplest way of 
proving a ``no loss of derivative'' result is to start from Cauchy's Lagrangian formulation of the
Euler equation \cite{Cau27, FV14}. This leads to a very simple and constructive proof of time-analyticity
of the Lagrangian trajectories \cite{FZ14, ZF14}.
Another recent related work \cite{CKV15}, uses a
Lagrangian formulation of the Euler equation to revisit some known results concerning 
the decay of the Eulerian analyticity radius in space variables 
\cite{Ben76, BB77, LO97, KV09, KV11-1, KV11-2, Zhe11} and to show the persistence 
of the spatial Lagrangian analyticity radius.

It is of interest to recall how the issue of Lagrangian regularity was tackled so far.
The temporal smoothness of Lagrangian particle paths, was first studied by Chemin (\cite{Che1}; see also \cite{Che2})
in Eulerian coordinates, using a paradifferential approach. He showed that Lagrangian trajectories
are $\mathscr{C}^\infty$--regular when initial conditions are  in the H\"older classes $\mathscr{C}^{\,m,\gamma}$ ($m\in\N-\{0\}$).
Chemin's proof relies on the estimate of the repeated action of the material derivative $\partial_t+v\cdot\nabla$ on
the pressure, which is expressed via a nonlinear pseudodifferential operator. In Refs. \cite{Ser92, Ser95-1, Ser95-2},
Serfati extended Chemin's results to the case of temporal analyticity and spatial H\"older regularity
($\mathscr{C}^{\,m,\gamma}$ with $m\in \N-\{0\}$). In the spirit of Refs., \cite{BB74} for an Euclidean space,
or \cite{EM70} for Riemannian manifolds, Serfati suggested to attack the problem as an abstract ordinary
differential equation (ODE) in a complex Banach space, holomorphic in time. The right-hand side of the ODE is
obtained by using the Biot--Savart-type singular integral representation of the pressure,
extended to an analytic functional on a complex Banach space. This requires
a careful analytic extension of the Newtonian potential singularity on a sector of the $n$-dimensional complex plane.
Serfati proved Banach--space analyticity of such a functional, and thus its Lipschitz property, by using the Cauchy integral
formula and showing boundedness of the functional in a complex Banach space.
More recently, Shnirelman  \cite{Shn12}
obtained another analytic abstract ODE formulation in complex Banach spaces, which makes use of
Arnold's geometric interpretation of the Euler equation \cite{Arn66}. The latter formulation uses the Biot--Savart-type singular
integral representation of the stream-function, and takes advantage of
Cauchy's vorticity formula in Lagrangian coordinates \cite{Cau27} to obtain a divergence-free representation of
the velocity field and discard the pressure. Serfati \cite{Ser95-2} and  Shnirelman \cite{Shn12} constructed
a time-analytic Lagrangian solution
to the ODE by using Picard successive approximations and a Banach fixed-point argument,
like in the proof of the Cauchy--Lipschitz--Picard theorem.
Another temporal analyticity result, which follows the spirit of Chemin's approach, is found in Ref. \cite{Gam94}.
Yet another approach, applicable to a class of hydrodynamic-like equations with nonlinear and nonlocal (collective)
interactions uses singular integral operators (SIO);
its proof makes crucial use of the  Fa\`a di Bruno formula \cite{FdB55, Gzy86, CS96, KP02} to deal with successive
time derivatives of the composition of the SIO kernel function with the flow map \cite{CVW15}.

All the previous regularity results were established for fluids without boundaries: either filling the whole space
or with periodic boundary conditions. In the presence of a boundary, using Eulerian coordinates, Kato
\cite{Kat00-1} showed $\mathscr{C}^\infty$--regularity in time of the Lagrangian trajectories
for H\"olderian classical solutions of the $n$-dimensional Euler 
equation in a smooth bounded domain with a  boundary of $\mathscr{C}^\infty$--regularity. Kato's
approach involved commutators of the material derivative and of the differential operators
arising from the Helmholtz--Hodge decomposition of the vector-valued Euler equation.
For these commutators, Kato obtained key estimates in both the interior of the domain and
--- even more important --- its boundary.
Of course these commutator estimates involve a representation of the boundary. Specifically 
the author used successive derivatives of the smooth function measuring
the distance between an interior point to the boundary. Nonlocal features, due to the pressure
term, were tackled by using classical elliptic regularity estimates. Refinements of Kato's
estimates have recently led to a proof of analyticity in time of the Lagrangian trajectories for
H\"olderian classical solution of the 3D Euler equation in the
presence of an analytic boundary \cite{GST12}.

The various methods we have briefly reviewed above are not constructive in the following sense:
they cannot be used to construct a high-order numerical scheme for integrating  the Euler equation
in either Eulerian or Lagrangian coordinates. A constructive approach follows from Cauchy's \cite{Cau27}
Lagrangian formulation. As shown in Refs. \cite{Cau27} this formulation leads to simple recursion relations among
the time-Taylor coefficients of the Lagrangian flow application mapping
initial fluid particle positions to current ones. 
It is then possible not only to prove time-analyticity, but also to construct efficiently very accurate
solutions to the Euler equation \cite{PZF16}. Here, we note that the 
terminology ``recursion relation'' is borrowed from the 
perturbation theory, as used in cosmology for over twenty years (see, e.g., \cite{BCGS02} and references
therein).

It is our main goal to show that such a constructive method can be extended to Euler flow with a boundary.
The outline of the paper is as follows. In Sec.~\ref{s:result}, we introduce our notation and define
our functional framework (Sec.~\ref{ss:def}). Then, in Sec.~\ref{ss:thm} we state our main result,
namely Theorem~\ref{thm:TUD}, a corollary of which is time-analyticity.
In Sec.~\ref{s:proof}, we present the proof of Theorem~\ref{thm:TUD} in two steps.
First, in Sec.~\ref{ss:const} we construct the time-Taylor coefficients by solving a sequence of non-homogeneous 
Dirichlet and Neumann boundary value problems. The boundary conditions on these are given by
boundary recursion relations stemming from the preservation of the boundary by the Euler flow (impenetrability).
Then, in Sec.~\ref{ss:conv}, we  deal with convergence issues and obtain bounds (e.g. in H\"older spaces) for time-Taylor coefficients,
from which the proof of Theorem~\ref{thm:TUD} follows.  
Sec.~\ref{s:concl} is devoted to concluding remarks and to discussing open problems.

\section{Smoothness of trajectories for an incompressible ideal fluid in a bounded domain}
\label{s:result}
In this section, we first define some notation and functional spaces. Then, we state our main theorem.
Actually, our functional framework consists of different classes of indefinitely differentiable functions,
that encompass {\em analytic}, 
{\em quasi-analytic}, and {\em non-quasi-analytic} classes.
We will perform our study in a general  ultradifferentiable class  \cite{Beu61, Rou62, LM72-3, Kom73} that we call the
``log-superlinear Fa\`a-di-Bruno'' class.

\subsection{Notation and functional framework}
\label{ss:def}
Let $N^\ast$ be the positive numbers of the set of natural integers $\N$. 
For a multi-index $\alpha=(\alpha_1,\ldots, \alpha_d)\in \N^d$ and  $x=(x_1,\ldots, x_d)\in \R^d$, we use the standard notation
$\alpha! := \alpha_1!\ldots \alpha_d!$, $\ \ |\alpha|:=\alpha_1+\ldots +\alpha_d$ and $x^\alpha:=x_1^{\alpha_1}\ldots x_d^{\alpha_d}$.
We use $\partial_i:=\partial/\partial x_i$, $\ \partial^\alpha:=\partial_1^{\alpha_1}\ldots \partial_d^{\alpha_d}$,
and we write $D^sf$ for the $s$-th Fr\'echet (functional) derivative of $f$.

We start with some definitions for ultradifferentiable classes, 
a generalization of Gevrey classes  \cite{Man42, LM72-3}.
Let $U$ be a domain in $\R^d$ and let $\mathfrak{B}$ be a Banach space endowed with the norm
$\| \cdot \|_\mathfrak{B}$. Let $f:U\longrightarrow \mathfrak{B}$ be an infinitely 
Fr\'echet-differentiable map. Then the $s$-th Fr\'echet derivative of $f$,
$D^sf:U\longrightarrow \mathcal{L}^s(U\times\ldots\times U, \mathfrak{B}) $,
is an element of $\mathcal{L}^s(U\times\ldots\times U, \mathfrak{B})$,
the set of $s$-linear bounded or continuous maps from $U\times\ldots\times U$ to $\mathfrak{B}$. The space
$\mathcal{L}^s(U\times\ldots\times U, \mathfrak{B})$ is endowed with the standard operator norm 
$\interleave \cdot \interleave:=\|\cdot\|_{\mathcal{L}^s(U\times\ldots\times U, \mathfrak{B})}$.

Let $\mathcal{M}:=\{ \mathcal{M}_s\}_{s\geq 0}$ be a sequence of positive numbers.
The ultradifferentiable class  $\mathcal{C}\{\mathcal{M}\} (U;\mathfrak{B})$ 
is defined  as the set of infinitely Fr\'echet-differentiable
functions $f:U\longrightarrow \mathfrak{B}$ such that  for any compact set $K\subset U$
there exist constants (depending on $f$) $R_f$, $C_f$ such that for all $s\in \N$, 
\begin{equation}
\label{UltraDif}
 \sup_{x\in K}\interleave D^sf(x) \interleave \leq C_fR_f^{-s}\mathcal{M}_s.
\end{equation}
The class $\mathcal{C}\{\mathcal{M}\}$ 
is invariant under multiplication by a constant, i.e.
$\mathcal{C}\{\lambda\mathcal{M}\}(U;\mathfrak{B})=\mathcal{C}\{\mathcal{M}\}(U;\mathfrak{B})$ for $\lambda >0$. 
There exists a fairly large set of ultradifferentiables functions \cite{Man42, Rou62, Kom73}.
Here we choose a class that we call the  ``log-superlinear Fa\`a-di-Bruno'' (LSL--FdB in short) class.
It involves a sequence of weights  $ {M}:=\{\mathcal{M}_s/s!\}_{s\geq 0}$ 
(and $M_0=\mathcal{M}_0$), in terms of which we have the following

\begin{definition} \label{def:LogSLFdB} The log-superlinear Fa\`a-di-Bruno class is the set of
functions satisfying (\ref{UltraDif}) with the following restriction on the weights: 
\begin{itemize}
\item [i)] differentiation stability:
\begin{equation}
\label{Diff-P}
\exists \  C_{{\rm d }}>0: \quad  M_{k+1} \leq  C_{{\rm d }}^k M_k, \quad  \forall k \in \N.
\end{equation}
\item [ii)] log-superlinearity: 
\begin{equation}
\label{log-SL-P}
M_k M_\ell \leq M_0 M_{k+\ell}, \quad \forall k, \, \ell \in \N.
\end{equation}
\item [iii)] (FdB)-stability: 
\begin{equation}
\label{FdB-P}
\exists \  0<C_{{\rm  FdB}}\leq 1 \ | \  \forall \alpha_i \in \N^\ast,\  \alpha_1+\ldots +\alpha_{\ell}=k:
\quad M_\ell M_{\alpha_1}\ldots  M_{\alpha_\ell} \leq  C_{{\rm FdB}}^kM_k.
\end{equation}
\end{itemize}
\end{definition}

A few remarks are now in order.

\begin{remark}
Using the Leibniz differentiation rules, log-superlinearity 
implies that the class $\mathcal{C}\{\mathcal{M}\} (U;\mathfrak{B})$  is an algebra with respect to 
pointwise multiplication. Using the Fa\`a di Bruno formula \cite{FdB55, Gzy86, CS96, KP02}, (FdB)-stability implies 
stability under composition in the class $\mathcal{C}\{\mathcal{M}\} (U;\mathfrak{B})$ 
(see the proof of Proposition 3.1 of Ref. \cite{RS14}). Finally, the
differentiability stability property implies closure under differentiation
in $\mathcal{C}\{\mathcal{M}\} (U;\mathfrak{B})$ \cite{Rou62, LM72-3, Kom73}. 
\end{remark}

\begin{remark}
Let us observe that this class is slightly larger than the log-convex class used in Ref. \cite{GST12}, 
in which the weights satisfy the following equivalent properties:
\begin{itemize}
\item[1.] $k\longmapsto \log M_k$ is convex.  
\item[2.] $M_k^2\leq M_{k-1}M_{k+1}, \ \forall k \in \N $.
\item[3.] The sequence  $\{ {M}_{s+1}/{M}_{s}\}_{s\geq 0}$  is nondecreasing.
\end{itemize}
It is straightforward to show that log-convexity implies log-superlinearity \eqref{log-SL-P} but the converse does not hold.
Moreover, it is shown in Lemma~2.2 of Ref. \cite{RS14} (see also \cite{RS16, KMR09}) that log-convexity implies (FdB)-stability
with an (FdB)-stability constant $C_{{\rm  FdB}}=\max\{1,\,M_1\}$. Since the class $\mathcal{C}\{\mathcal{M}\}$ 
is invariant under multiplication by a constant,
we can normalize the sequence ${M}=\{M_k\}_{k\geq 0}$ by an arbitrary positive constant.
Therefore, dividing ${M}$ by $M_1$  implies  $C_{{\rm  FdB}}=\max\{1,\,M_1\}=1$. In other words
log-convexity implies (FdB)-stability \eqref{FdB-P}.
Let us note that usually in the literature \cite{Beu61, Rou62, LM72-3, Kom73, Thi08, KMR09, KMR11, KMR15, RS14, RS16},
ultradifferentiable classes assume log-convexity in order to achieve stability under composition.
We refer the reader to \cite{RS14, RS16} and references therein for a discussion of the role of
log-convexity in achieving stability under composition.
\end{remark}

\begin{remark}
Some well-known classes of functions belong to the LSL--FdB class. 
The first one is the analytic functions class which corresponds to $M_k=1$ (e.g. \cite{Man42, Car61}). The second one is
the Gevrey class (see, e.g., \cite{Man42, LM72-3} and references therein)
which corresponds to $M_k=(k!)^r$ with $r>0$. In the usual technical sense, the former
is also quasi-analytic and the latter is not \cite{Rud87, Beu61, Hor83, Thi08, KMR09, KMR11} (see also \cite{KMR15}).  
Recall that ultradifferentiable classes are divided into quasi-analytic and non-quasi-analytic classes,
depending on whether the map to infinite Taylor expansions is injective on the class or not (see, e.g., \cite{Thi08, KMR09, KMR11}).
Denjoy \cite{Den21} and Carleman \cite{Car26} have given some criteria for distinguishing
quasi-analytic classes from non-quasi-analytic ones. Typically, the Denjoy--Carleman theorem \cite{Den21, Car26}
states that  $\mathcal{C}\{\mathcal{M}\} (U;\mathfrak{B})$ lacks quasi-analyticity if and only if $\ \sum_{k\geq 0}M_k/M_{k+1} <\infty$
(see, e.g., \cite{Rud87, Hor83, Kom73, Thi08, KMR09, KMR11} for contemporary proofs).
We refer the reader to \cite{Beu61, Rou62,  Kom73,  Hor83, Rud87, Thi08, KMR09, KMR11, KMR15} and references therein 
for establishing various equivalent Denjoy--Carleman-type criteria.
\end{remark}

The ultradifferentiable framework will be used here for studying the regularity in time. 
It remains to choose
a functional space to measure the regularity in space of the initial data, i.e. the regularity with respect
to the Lagrangian spatial variable. Two requirements are needed for this functional space. First, as we shall see, it
must be an algebra with respect to the pointwise multiplication. Second, it must be compatible
with local-in-time well-posedness results of the 3D incompressible Euler equation (see Theorem \ref{thm:K2000} below). 
For convenience, we choose H\"older spaces, but our proof can be extended, for instance, to Sobolev and Besov spaces 
(with suitable indexes of regularity), or ultradifferentiable spaces (e.g.,  analytic functions). 
We recall the definition of the H\"older spaces $\mathscr{C}^{\,m,\gamma}$
used here and 
an associated basic property. Let $\Omega$ be
a bounded domain in $\R^d$ with $d\geq 1$. Let $\mathscr{C}^{\,m}(\overline{\Omega})$ be the space of $m$ times 
continuously differential functions on the closure $\overline{\Omega}$, equipped with the norm
\begin{equation*}
\label{norm-C-m}
\|f\|_m \equiv \|f\|_{\mathscr{C}^{\,m}(\overline{\Omega})} := \max_{0\leq |\alpha|\leq m}\  \sup_{x\in \overline{\Omega}}|\partial^\alpha f(x)|. 
\end{equation*}
For $0<\gamma<1$ and $m\in \N$, the H\"older space  $\mathscr{C}^{\,m,\gamma}(\overline{\Omega})$ is defined as
the subspace of $\mathscr{C}^{\,m}(\overline{\Omega})$ consisting of those functions $f$ for 
which, for $0<|\alpha|\leq m$, $\,\, \partial^\alpha f$ satisfies in $\Omega$ a H\"older condition of exponent
$\gamma$. The space  $\mathscr{C}^{\,m,\gamma}(\overline{\Omega})$ is endowed with the norm
\begin{equation*}
\label{norm-Holder-m}
\|f\|_{m,\gamma }\equiv \|f\|_{\mathscr{C}^{\,m,\gamma}(\overline{\Omega})} := \|f\|_m +
\max_{0\leq |\alpha|\leq m} \  \sup_{\substack{x,y \in \overline{\Omega} \\ x\neq y}} \frac{|\partial^\alpha f(x) - \partial^\alpha f(y) |}{|x-y|^\gamma}. 
\end{equation*}
The space  $\mathscr{C}^{\,m,\gamma}(\overline{\Omega})$ is an algebra with respect to  the pointwise multiplication 
(see, e.g., Chapter 4 of \cite{MB02}, or Chapter 4 of \cite{GT98}), 
i.e. there exists a constant $C_a:=C_a(m)$, which depends on $m$, but is independent of $\gamma$, such that
\begin{equation}
\label{algebra-Holder-m}
\|fg\|_{m,\gamma } \leq C_a \|f\|_{m,\gamma } \|g\|_{m,\gamma }, \quad \forall f, \, g \in  \mathscr{C}^{\,m,\gamma}(\overline{\Omega}).
\end{equation}
Finally we introduce the notation $\mathcal{D}'(\Omega)$ for the space of (Schwartz) distribution on $\Omega$.
The space $\mathcal{D}'(\Omega)$ is defined in a standard way as the topological dual of the space of
indefinitely differentiable functions with compact support on $\Omega$.
\subsection{Main theorem}
\label{ss:thm}
Before tackling issues of Lagrangian time analyticity and variants thereof, we state a key result of Kato
about the  initial-value Euler problem in the presence of a boundary  \cite{Kat00-1}. Kato's result
applies in particular to an incompressible ideal fluid filling a bounded
simply connected regular domain $\Omega \subset \R^3$, as we shall assume throughout 
this paper. The relevant equations (so far in Eulerian coordinates) are: 
\begin{alignat}{4}
\label{euler-1} 
\partial_t v + (v\cdot \nabla) v + \nabla p & = & 0, &\ \quad x\in\Omega, &\ \  t\in ]-T,T[ \\
\label{euler-2} 
\nabla\cdot u &= & 0,  &\ \quad x\in\Omega, &\ \ t\in]-T,T[ \\
\label{euler-3} 
v_{|_{t=0}}&= & \ v_0, &\ \quad x\in\Omega &\\
\label{euler-4} 
v\cdot \nu &= &0,  &\ \quad x\in\partial \Omega, &\,\ \  t\in]-T,T[.
\end{alignat}
Here $v$ is the velocity,  $p$ the pressure and $\nu$ the outer normal at the 
impermeable boundary $\partial\Omega$.

We then have the following 
\begin{theorem}(Kato  \cite{Kat00-1}). 
\label{thm:K2000}
Suppose that $(v,p)$ is a solution of \eqref{euler-1}--\eqref{euler-4}, where 
$v\in  \mathscr{C}\left(]0,T[;\mathscr{C}^{\,1,\gamma}(\overline{\Omega})\right)$ with $\gamma\in (0,1)$, and
$p\in \mathcal{D}'(]0,T[\times\Omega)$. Then we have, for all $s\in \N^\ast$, 
\begin{equation}
\label{K2000_est_1}
(\partial_t + v\cdot \nabla)^{s}v\,\in \mathscr{C}\left(]0,T[;\mathscr{C}^{\,1,\gamma}(\overline{\Omega})\right),
\quad \nabla(\partial_t + v\cdot \nabla)^{s-1}p\,\in \mathscr{C}\left(]0,T[;\mathscr{C}^{\,0,\gamma}(\overline{\Omega})\right),
\end{equation}
and
\begin{equation}
\label{K2000_est_2}
\|(\partial_t + v\cdot \nabla)^{s}v(t)\|_{1,\gamma}\ +\ 
\|\nabla(\partial_t + v\cdot \nabla)^{s-1}p(t)\|_{0,\gamma} \leq C_s \|v(t) \|_{1,\gamma}^{s+1},
\quad \forall t\in]0,T[,
\end{equation}
where $C_s$ are constants depending on $\gamma$,  $\, s$ and $\Omega$.
\end{theorem}

\begin{remark}
\label{rem:LWP}
We choose here the H\"older space $\mathscr{C}^{\,1,\gamma}(\overline{\Omega})$ but, of course, Kato's theorem
also holds for H\"older spaces $\mathscr{C}^{\,m,\gamma}(\overline{\Omega})$ with $m\geq 1$ (see \cite{Kat00-1}).
Theorem  \ref{thm:K2000} relies on existence, uniqueness and regularity properties of
the local-in-time classical solution of the Euler equations \eqref{euler-1}--\eqref{euler-4} in
H\"older spaces (see, e.g., Theorem I of \cite{Kat00-2}) with initial condition 
$v_0\in \mathscr{C}^{\,1,\gamma}(\overline{\Omega})$ and such that $T=C_\ast(\Omega)/\|v_0\|_{\mathscr{C}^{\,1,\gamma}(\overline{\Omega})}$.
Moreover, Theorem \ref{thm:K2000} can be extended for,
say, initial conditions in any Sobolev space $H^s(\Omega)$ with $s>5/2$ and any inhomogeneous Besov space
$B_{p,q}^s(\Omega)$, with $1\leq p,\, q\leq \infty$ and with $s>3/p +1$  or
$s\geq 3/p+1$ if $q=1$, so that $B_{p,q}^s(\Omega)$ is continuously
embedded in the Lipschitz space ${\rm Lip}(\Omega)$
(see \cite{Dut03}).
Let us recall that it is still not known whether the classical solutions of Theorem \ref{thm:K2000}
remain smooth for all times or ``blow-up'' in finite time. Let us also point out the recent work
by Bardos and Titi \cite{BT10} showing that the 3D Euler equations are not well-posed in the
H\"older spaces $\mathscr{C}^{\,0,\gamma}(\overline{\Omega})$ for $\gamma\in (0,1)$.
\end{remark}

For the statement of our main theorem we need to introduce the Lagrangian map $X$
defined on $]-T,T[\times \Omega$ by 
\begin{equation}
\label{LagMap}
\partial_t X(t,a) = v(t,X(t,a)) \ \ \ {\rm and } \ \ \ X(0,a)=a.
\end{equation}
The Lagrangian map $X(t,\cdot):\Omega \longrightarrow \Omega $ can been seen as a continuous one-parameter
group of
volume and orientation preserving diffeomorphisms defined on $\Omega$ \cite{Arn66}. Fields
expressed in terms of the Lagrangian label $a$ and time $t$ are said to be given in Lagrangian
coordinates. The Lagrangian gradient, with respect to $a$, is denoted $\nabla^{\rm L}$ and its components
$\partial_i^{\rm L}$. 

The main result of this paper states roughly that the smoothness of the Lagrangian characteristic curves is only
limited by the smoothness of the domain boundary. More precisely, we have
\begin{theorem}
\label{thm:TUD}
Assume that the hypotheses of  Theorem \ref{thm:K2000} hold, and in addition that the boundary $\partial \Omega$
belongs to $\mathcal{C}\{\mathcal{M}\}$, where $\mathcal{M}:=\{ s!M_s\}_{s\geq 0}$, with the sequence 
$\{M_s\}_{s\geq 0}$ satisfying Definition \ref{def:LogSLFdB} (log-superlinear Fa\`a-di-Bruno class).
Then there exists a time  $T=C_\ast(\Omega)/\|v_0\|_{\mathscr{C}^{\,1,\gamma}(\overline{\Omega})}$ such that
the Lagrangian  map $X$ satisfies
\[
X \in \mathcal{C}\{\mathcal{M}\}\left(]-T,T[; \mathscr{C}^{\,1,\gamma}(\overline{\Omega})\right). 
\]
\end{theorem}

From  Theorem \ref{thm:TUD} we infer, by specializing to 
${\mathcal M} = \{s!\}_{s\geq 0}$ or, equivalently to $M_s =1, \, \forall s\geq 0$, the following analyticity\\[-3ex]
\begin{corollary}
\label{cor:anal}
Assume that the hypotheses of  Theorem \ref{thm:TUD} hold.  If the boundary $\partial \Omega$
is analytic (resp. Gevrey of order $r>0$), then  Lagrangian map $X$ is analytic
(resp. Gevrey of order $r>0$)  from $]-T,T[$ to $\mathscr{C}^{\,1,\gamma}(\overline{\Omega})$.
For the analytic case, the  Lagrangian map $X$ admits a convergent time-Taylor expansion
around time $t=0$.
\end{corollary}

\section{Proof of Theorem \ref{thm:TUD}}
\label{s:proof}
Here, we give a proof of Theorem \ref{thm:TUD},  divided into two steps. Sec.~\ref{ss:const}
is devoted to a recursive construction of a solution, for the Lagrangian map, as a formal 
time-Taylor  expansion, without specifying any functional setting.  
The second step (Sec.~\ref{ss:conv})  is dedicated to convergence issues  of such formal 
time-expansions within the functional framework defined in Sec.~\ref{ss:def}.

\subsection{Construction of a solution as a formal time series}
\label{ss:const}
Our starting point is not the usual Eulerian formulation of the Euler
equations (\ref{euler-1})--(\ref{euler-4}) for incompressible ideal fluid, but
a little-known  Lagrangian formulation due to Cauchy \cite{Cau27}, together
with the incompressibility condition and the preservation by the Lagrangian
flow of a prescribed rigid boundary $\partial \Omega$, represented 
by the equation $\mathcal{S}(x)=0$ ($\mathcal{S}:\R^3\rightarrow \R$
being a 
map of appropriate regularity). The relevant equations are:
\begin{alignat}{6}
\sum_{k=1}^3\nabla^{\rm L}\dot{X}_k(t,a)\times \nabla^{\rm L}{X}_k(t,a)& = &&\ \omega_0(a), 
 && \quad \forall a\in \Omega,\label{Cau} \\[-1.4ex]
\det (\nabla^{\rm L}{X}(t,a)) &=&&\ 1, && \quad \forall a\in \Omega, \label{Lag}\\[0.5ex]
\mathcal{S}(X(t,a))& = &&\ 0, && \quad \forall a\in \partial \Omega.\label{InvSurf}
\end{alignat}
Here, the Lagrangian variable $X$ (resp. $a$) stands for the current (resp. initial) Lagrangian position of 
fluid particules. The initial vorticity $\omega_0$ is defined as usual, i.e. $\omega_0= \nabla^{\rm L}\times v_0$,
where $\nabla^{\rm L}$ denotes the gradient in the Lagrangian variables. The dot denotes the Lagrangian
time derivative, while $\nabla^{\rm L}{X}$ stands for the Jacobian matrix with entries $\partial_i^{\rm L}{X}_j$.
The left-hand-side of (\ref{Cau}), which is obviously time-invariant, is usually
referred to as the ``Cauchy invariants.''

Now, we observe that the assumptions of our Theorem~\ref{thm:TUD} are always
stronger than those of Kato's Theorem~\ref{thm:K2000}, given in Section~\ref{ss:thm}. Thus, the initial velocity
is at least $\mathscr{C}^{1,\gamma}$ and the boundary $\partial\Omega$ 
has at least $\mathscr{C}^\infty$ regularity, from which
follows that the Lagrangian map is $\mathscr{C}^\infty$
in time.  The recursion relations among time-Taylor coefficients derived
below could be obtained from the basic Lagrangian equation
\eqref{Cau}--\eqref{InvSurf} by successive time-differentiations and
use of standard relations such as the Fa\`a di Bruno formula for
differentiation of compositions. This procedure can be somewhat algebraized
by using formal time-Taylor series without worrying about convergence issues.
Indeed, as we shall see, the relations we will obtain always involve a finite
number of terms. Of course, once we address issues such as analyticity,
studying the convergence of such formal series becomes important. This
will however be postponed to Sec.~\ref{ss:conv}, where we will derive various
functional estimates. 

For the implementation, it is convenient, following Refs.~\cite{FZ14, ZF14}, 
to introduce the displacement field $\xi=\xi(t,a)=X-a$, in terms of which,
\eqref{Cau}--\eqref{Lag} become 
\begin{eqnarray}
&& \nabla^{\rm L} \times \dot{\xi} + \sum_{1\leq k\leq 3} \nabla^{\rm L}  
\dot{\xi}_k \times  \nabla^{\rm L} {\xi}_k= \omega_0, \label{CauDep}\\
&&  \nabla^{\rm L} \cdot \xi +  \sum_{1\leq i<j \leq 3}\left(  \partial_i^{\rm L} {\xi}_i \partial_j^{\rm L} {\xi}_j
- \partial_i^{\rm L} {\xi}_j \partial_j^{\rm L} {\xi}_i \right)  
+ \det  (\nabla^{\rm L}{\xi}) =0.\label{LagDep}
\end{eqnarray} 
We now use the formal Taylor series for the displacement
\begin{equation}
\label{TTSExi}
\xi(t,a)=\sum_{s> 0}\xi^{(s)}(a)t^s.
\end{equation}
Substituting \eqref{TTSExi} into \eqref{CauDep}--\eqref{LagDep} and collecting
terms of the same power $s>0$, we obtain, after a symmetrization of the sums
\begin{eqnarray}
  \nabla^{\rm L} \times {\xi}^{(s)} &=&\omega_0\delta_{1s} -\frac{1}{2}
\sum_{\substack{ 1\leq k\leq 3\\ 0<m<s}} \frac{2m-s}{s}  \nabla^{\rm L} {\xi}_k^{(m)} \times\nabla^{\rm L} {\xi}_k^{(s-m)} \label{rotxis}  
\\
\nabla^{\rm L} \cdot \xi^{(s)} &=& \sum_{\substack{ 1\leq i<j\leq 3\\ 0<m<s}} 
\left(\partial_i^{\rm L} {\xi}_j^{(m)} \partial_j^{\rm L} {\xi}_i^{(s-m)} 
-\partial_i^{\rm L} {\xi}_i^{(m)} \partial_j^{\rm L} {\xi}_j^{(s-m)}\right)  \nonumber\\
&&
-\frac{1}{6} \sum_{l+m+n=s} \varepsilon_{ijk}\varepsilon_{abc} \partial_i^{\rm L} {\xi}_a^{(l)}
 \partial_j^{\rm L} {\xi}_b^{(m)}  \partial_k^{\rm L} {\xi}_c^{(n)}. \label{divxis} 
\end{eqnarray}
Here, $\varepsilon_{ijk}$ stands for the unit antisymmetric tensor and $\delta_{ij}$ is the Kronecker symbol.
In \eqref{divxis}, repeated indices are implicity summed over.

We observe that \eqref{rotxis}--\eqref{divxis}, which prescribe the
Lagrangian curl and
divergence of ${\xi}^{(s)}$,  constitute a Helmholtz--Hodge problem (see, e.g.,
\cite{FT78, GR86, ABDG98}). This is the point where it is
essential to handle the boundary conditions, stemming from \eqref{InvSurf}.
Indeed, the 3D Helmholtz--Hodge decomposition of the vector field $\xi^{(s)}$ 
consists in expressing it as follows: 
\begin{equation}
\label{HH_exp}
\xi^{(s)} = \nabla^{\rm L} \varphi^{(s)} \,+\, \nabla^{\rm L} \times  \Phi^{(s)},
\quad \ {\rm in }\ \,\overline{\Omega},
\end{equation}
where $\varphi^{(s)}$ (resp.  $\Phi^{(s)}$) is the scalar (resp. vector) potential
(see, e.g., \cite{FT78, GR86, ABDG98}). Without loss of generality one can assume
the gauge condition
\begin{equation}
  \label{gauge}
  \nabla^{\rm L} \cdot  \Phi^{(s)}=0, \quad  \ {\rm in }\ \,\overline{\Omega}.
\end{equation}  
By taking the divergence and the curl of (\ref{HH_exp}) and using \eqref{gauge},
it is found that the scalar potential
$\varphi^{(s)}$ and the vector potential $\Phi^{(s)}$ satisfy, in the interior
of $\Omega$, the following non-homogeneous Poisson equations
\begin{equation}
\label{Laplace_phiPhi}
\Delta^{\rm L} \varphi^{(s)}= \nabla^{\rm L} \cdot \xi^{(s)}, \quad  \ \
\Delta^{\rm L} \Phi^{(s)} = -\,\nabla^{\rm L} \times {\xi}^{(s)}, \quad \ {\rm in} \ \,\Omega,
\end{equation}
where $\Delta^{\rm L}$ is the Lagrangian Laplacian. An important observation is that the boundary traces of  the scalar potential
$\varphi^{(s)}$ and of the vector potential $\Phi^{(s)}$ on $\partial \Omega$ are obtainable from
the flow-invariance of the boundary, expressed by the constraint \eqref{InvSurf}.
Indeed, for a given $a\in \partial \Omega$, let us define
\begin{equation}
\label{def_f}
f(t)=\mathcal{S}(X(t,a)), 
\end{equation}
which has to vanish, since $X(t,a)$ remains on $\partial \Omega $. We now expand \eqref{def_f}
in a (possibly formal) time-Taylor series,
\begin{equation}
\label{SurfTayl}
f(t)=\sum_{s\geq 0} f^{(s)}(0) \frac{t^s}{s!}, \quad \quad \forall a\in \partial \Omega,
\end{equation}
where the $s$-th derivative $f^{(s)}(0)$ can be expressed by using the Fa\`a di Bruno formula \cite{FdB55, Gzy86, CS96}:
\begin{equation}
\label{TaylCoefFdB}
 f^{(s)}(0) = s! \sum_{1\leq |\beta| \leq s} \partial^\beta\mathcal{S}(a)
\sum_{i=1}^s \sum_{P_i(s,\beta)} \prod_{j=1}^i 
\frac{(\xi_1^{(\ell_j)})^{k_{j1}}}{k_{j1}!}
\frac{(\xi_2^{(\ell_j)})^{k_{j2}}}{k_{j2}!}
\frac{(\xi_3^{(\ell_j)})^{k_{j3}}}{k_{j3}!},
\end{equation}
for $s>0$ and $f(0)=\mathcal{S}(a)=0$.
In (\ref{TaylCoefFdB}) the set $P_i(s,\beta)$ is given by
\begin{multline}
\label{PartitionSet} 
P_i(s,\beta)= \Bigg\{
(\ell_1,\ldots,\ell_i), \ (k_1,\ldots,k_i); \ \ 0<\ell_1 < \ldots < \ell_i;  \\
\ \ |k_j|>0, \ j\in [1,i];\ \
\sum_{j=1}^ik_j=\beta, \ \ \sum_{j=1}^i|k_j|\ell_j=s
\Bigg\}.
\end{multline}
From  (\ref{SurfTayl})--(\ref{TaylCoefFdB}) we infer the relations
\begin{equation}
\label{NBC_1}
 \xi^{(s)}(a)\cdot \nabla \mathcal{S}(a) = - \sum_{1<|\beta| \leq s} \partial^\beta\mathcal{S}(a)
\sum_{i=1}^s \sum_{P_i(s,\beta)} \prod_{j=1}^i 
\frac{(\xi_1^{(\ell_j)})^{k_{j1}}}{k_{j1}!}
\frac{(\xi_2^{(\ell_j)})^{k_{j2}}}{k_{j2}!}
\frac{(\xi_3^{(\ell_j)})^{k_{j3}}}{k_{j3}!}, \quad \forall q\in\partial\Omega,
\end{equation}
which play a key role in obtaining recursion relations for the Taylor coefficients
$\xi^{(s)}$. Let us already observe that the right-hand side of \eqref{NBC_1}
involves only Taylor coefficients of order less than $s$. This follows from
the definition \eqref{PartitionSet} of the multi-index set $P_i(s,\beta)$.

Using  \eqref{NBC_1} and the Helmholtz--Hodge decomposition \eqref{HH_exp}--\eqref{Laplace_phiPhi},
we shall now show that the determination of the Taylor coefficient $\xi^{(s)}$,
knowing all those of lower order, reduces to a vector-valued non-homogeneous
Dirichlet problem and to a non-homogeneous Neumann problem.
Let us first note that $\nu$, the normal vector to the boundary
$\partial \Omega$, may be expressed in terms of the defining-boundary map $\mathcal{S}$ as
\begin{equation}
\label{def_nu_by_S}
\nu = \frac{\nabla \mathcal{S}} {|  \nabla \mathcal{S}|_2},
\end{equation}
where $|\cdot|_2$ is the Euclidean norm and $|\nabla \mathcal{S}|_2>0$.
Thus the left-hand side of \eqref{NBC_1} is expressed in terms of the normal
component of the Taylor coefficient of order $s$. Upon using \eqref{HH_exp},
we obtain a surface condition involving both the scalar  potential
$\varphi^{(s)}$ and the vector potential $\Phi^{(s)}$. We choose to impose
the vanishing of the vector potential $\Phi^{(s)}$ on $\partial \Omega$;
this implies the vanishing of $  \nu \cdot \nabla^{\rm L} \times \Phi^{(s)}$.
In this way the equations for the two potential decouple and take the following form.
The scalar potential $\varphi^{(s)}$ satisfies the
non-homogeneous Neumann boundary value problem
\begin{equation}
\label{NeumannBVP}
\left \{
\begin{tabular}{llll}
$\Delta^{\rm L} \varphi^{(s)}$ &$=$&  $\nabla^{\rm L} \cdot {\xi}^{(s)} $ & {\rm in} $\ \Omega$ \vspace{3pt}\\ 
$\displaystyle{\partial_{\nu}^{\rm L}  \varphi^{(s)}}$&=&$\xi^{(s)}\cdot \nu $ & {\rm on} $\ \partial \Omega$,
\end{tabular}
\right.
\end{equation}
where it is understood that the right-hand sides of (\ref{NeumannBVP}) are taken from
(\ref{divxis}) and  (\ref{NBC_1})--(\ref{def_nu_by_S}), and thus involve only Taylor
coefficients $\{\xi^{(s)} \}_{0\leq s \leq s-1}$. As to
the vector potential $\Phi^{(s)}$, it satisfies the
non-homogeneous Dirichlet boundary value problem
\begin{equation}
\label{DirichletBVP}
\left \{
\begin{tabular}{llll}
$\Delta^{\rm L} \Phi^{(s)}$ &$=$&  $-\,\nabla^{\rm L} \times {\xi}^{(s)}$ &{\rm in} $\ \Omega$ \vspace{3pt} \\ 
$\Phi^{(s)}$&=&$0$ &{\rm on} $\ \partial \Omega$,
\end{tabular}
\right.
\end{equation}
where the right-hand side of (\ref{DirichletBVP}) is given by (\ref{rotxis}).

It is straightforward to solve
the  Helmholtz--Hodge decomposition \eqref{HH_exp}
at the order $s=1$. Indeed, from  (\ref{rotxis}) we obtain that
$\nabla^{\rm L} \times (\xi^{(1)}-v_0)=0$, which implies that there exists a function
$f$ such that $\nabla^{\rm L} f=\xi^{(1)}-v_0$. Then,  from  (\ref{divxis})
we obtain that $\nabla^{\rm L} \cdot \xi^{(1)}=0$, which together with the initial
incompressibility condition $\nabla \cdot v_0=0$ implies $\Delta^{\rm L} f= \nabla^{\rm L} \cdot (\xi^{(1)}-v_0)=0$.
Adding the consistent homogeneous Neumann condition $\nu \cdot \nabla^{\rm L} f =0$ on $\partial \Omega$,
to the homogeneous Laplace equation $\Delta^{\rm L} f=0$ in $\Omega$, we obtain $f=0$ in $\overline{\Omega}$,
i.e.
\begin{equation}
\label{recur1}
\xi^{(1)}=v_0.
\end{equation}
This  is the starting point of the recursive procedure 
for the Taylor coefficients.  

A few remarks are now in order.

\begin{remark}
\label{rem:uniq_DNBVP}
It is standard to prove existence and uniqueness (sometimes up to an additive constant)
of the non-homogeneous Dirichlet and Neumann boundary value problems, when
the right-hand side and the boundary terms are in suitable functional spaces. Such
issues are postponed to Sec.~\ref{ss:conv}, where we use such spaces for obtaining
a priori estimates.
\end{remark}

\begin{remark}
\label{rem:uniq_xis}
We have assumed, for convenience that $\Phi^{(s)}=0$ at the boundary. There are
other possible choices (see, e.g., \cite{GR86, ABDG98}), but it is not difficult to
show from \eqref{rotxis}--\eqref{divxis} and \eqref{NBC_1} that all lead to the same
solution $\xi^{s}$.
\end{remark}

\begin{remark}
\label{rem:uniq_D}
We have used a representation of the boundary $\partial \Omega$
by the vanishing of some map $\mathcal{S}$. This is definitely
not unique. Furthermore the map $\mathcal{S}$ may be defined either only on the
neighborhood of $\partial \Omega$ or more globaly \cite{BB74, Kat00-1}.
Will these choices lead to the same time-Taylor coefficients ? This can
be proved indirectly using Kato's Eulerian classical solution \cite{Kat00-2}, but it would be useful to have a direct Lagrangian
proof.
\end{remark}

\subsection{Convergence analysis of the formal time series}
\label{ss:conv}
Here, we prove that the Lagrangian map $X$ is time-ultradifferentiable in the log-superlinear
Fa\`a di Bruno class
$\ \mathcal{C}\{\mathcal{M}\} \left(]-T,T[;\mathscr{C}^{\,1,\gamma}(\overline{\Omega})\right)$,
defined in Sec.~\ref{ss:def}. As a by-product, 
we obtain the convergence of the time-Taylor expansion (\ref{TTSExi}) of the Lagrangian map $X$,
which was explicitly constructed in Sec.~\ref{ss:const}.

Let us introduce the generating function
\begin{equation}
\label{GenFunc}
\zeta(t) = \sum_{s> 0} \|\xi^{(s)}\|_{1,\gamma}\,M_s^{-1}\, t^s.
\end{equation}
Then, showing that $X\in \mathcal{C}\{\mathcal{M}\} \left(]-T,T[;\mathscr{C}^{\,1,\gamma}(\overline{\Omega})\right)$
for some $T>0$ is equivalent to proving that $\zeta(t)$ is uniformly bounded  on $]-T,T[$, or that 
the right-hand side of (\ref{GenFunc}) converges uniformly on $]-T,T[$. 
To simplify the exposition, we first prove Theorem \ref{thm:TUD} for initial data that are analytic
in the space variable but, nevertheless, as we shall see, we will obtain a priori estimates
that depend only on $\| v_0\|_{1,\gamma}$.
This is in fact enough for obtaining the general result, since by standard approximation methods, we can regularize $v_0$ to render the initial
data analytic in space. Then, from  results of Refs. \cite{Ben76, BB77, Del85}, we obtain
a unique solution that is analytic in time and space for some time.
In particular, the displacement $\xi$ is analytic both in time and
space and thus the time-Taylor series (\ref{GenFunc}) converges in a non-empty disk.
Then, it will be clear from the proof below  that the estimates are uniform
for the regularized solution. We thus can pass to the limit,
keeping the desired a priori estimates by standard arguments, eventually being
able to handle initial data $v_0\in \mathscr{C}^{1,\gamma}(\overline{\Omega})$.

Let us now derive the a priori estimates. For this we will use Schauder's regularity estimates in H\"older spaces
for the non-homogeneous Neumann and Dirichlet boundary
 value problems, needed to solve the Helmholtz--Hodge decomposition
 problems that
appear in the previous section. Specifically, for the Dirichlet problem we shall use the
H\"older estimates found in Refs. \cite{LU68, Mir70, GT98}. As to the more delicate
Neumann problem, there are results using Sobolev spaces in Refs. \cite{LU68, LM72-1, GT98},
but for the H\"older case we shall use the recent result of Nardi \cite{Nar14}.

We now wish to estimate, in the space $\mathcal{C}^{1,\gamma}$, the Taylor coefficient $\xi^{(s)}$
of order $s$, given by the Helmholtz--Hodge decomposition \eqref{HH_exp}. For this, we need
$\mathcal{C}^{2,\gamma}$ estimates on the potentials $\varphi^{(s)}$ and  $\Phi^{(s)}$, solutions
of the Neumann problem \eqref{NeumannBVP} and the Dirichlet problem \eqref{DirichletBVP}. 
The aforementioned estimates imply that there exist two constants $C_{\rm D}=C_{\rm D}(\Omega,\gamma)$ and 
$C_{\rm N}=C_{\rm N}(\Omega,\gamma)$ such that 
\begin{eqnarray}
  \|\xi^{(s)}\|_{1,\gamma} &\leq & \left\|\nabla^{\rm L} \varphi^{(s)}\right\|_{1,\gamma} +
  \left\| \nabla^{\rm L}\times \Phi^{(s)}\right\|_{1,\gamma} \nonumber\\[0.9ex]
&\leq& C_{\rm num}\left(\left\| \varphi^{(s)}\right\|_{2,\gamma} + \left\|\Phi^{(s)}\right\|_{2,\gamma}\right) \nonumber\\[0.9ex]
  &\leq& C_{\rm DN} \left(\left \|\nabla^{\rm L} \cdot \xi^{(s)}\right\|_{0,\gamma} +
  \left\|\nabla^{\rm L} \times \xi^{(s)}\right\|_{0,\gamma} + \left\|\xi^{(s)}\cdot \nu\right\|_{1,\gamma}
\right), \label{xisest_1}
\end{eqnarray}
where $ C_{\rm DN}:=C_{\rm num}\max\{  C_{\rm D},  C_{\rm N}\}$ and $C_{\rm num}$  is a pure numerical constant. Let us note that 
the operators $\, \partial_i^{\rm L}\partial_j^{\rm L}\Delta_{\rm L}^{-1}\,$ and  $\, \varepsilon_{ik\ell}
\partial_j^{\rm L}\partial_k^{\rm L}\Delta_{\rm L}^{-1}(\cdot)_\ell \,$
are Calder\'on-Zygmund operators of degre zero, which are  continuous endomorphisms in the
homogeneous H\"older space $\mathcal{C}^{0,\gamma}(\Omega)$, with a $\gamma$-depending continuity constant
varying as $1/\gamma$ (see, e.g., Chapter 4 of \cite{MB02}, or Chapter 4 of \cite{GT98}). Therefore, the constants 
$C_{\rm D}$ and $C_{\rm N}$, and consequently
$ C_{\rm DN}$ vary in $\gamma$ as $1/\gamma$.  
We then  must estimate the right-hand side of (\ref{xisest_1}).
Using the algebra property (\ref{algebra-Holder-m}) and
the superlinearity property (\ref{log-SL-P}), we find that 
the divergence of $\xi^{(s)}$, given by the right-hand side of (\ref{divxis}), has the following bound
\begin{equation*}
\label{divxis_est_1}  
\frac{\left\|\nabla^{\rm L} \cdot \xi^{(s)}\right\|_{0,\gamma}}{M_s} 
\leq 6C_{a}M_0 \sum_{0<m<s}  \frac{\| \xi^{(m)}\|_{1,\gamma}}{M_m}  \frac{\| \xi^{(s-m)}\|_{1,\gamma}}{M_{s-m}}
+  6C_{a}^2M_0^2 \sum_{l+m+n=s} \frac{\| \xi^{(l)}\|_{1,\gamma}}{M_l} \frac{\| \xi^{(m)}\|_{1,\gamma}}{M_m}\frac{\| \xi^{(n)}\|_{1,\gamma}}{M_n},
\end{equation*}
where $C_a$ is appearing in the algebra bound (\ref{algebra-Holder-m}).
It follows that
\begin{equation}
\label{divxis_est_2}  
\sum_{s>0}\left\|\nabla^{\rm L} \cdot \xi^{(s)}\right\|_{0,\gamma}\, M_s^{-1}\, t^s
\leq 6C_{a} M_0\left(\zeta^2 + C_aM_0\zeta^3 \right).
\end{equation}
In an exactly similar way we can bound the curl of $\xi^{(s)}$, given by the right-hand side of (\ref{rotxis}): 
\begin{equation*}
\label{rotxis_est_1}  
\frac{\left\|\nabla^{\rm L} \times \xi^{(s)}\right\|_{0,\gamma}}{M_s} 
\leq \frac{1}{M_1}\|\omega_0\|_{0,\gamma}\delta_{1s} +
\frac{3}{2}C_{a}M_0 \sum_{0<m<s}  \frac{\| \xi^{(m)}\|_{1,\gamma}}{M_m}  \frac{\| \xi^{(s-m)}\|_{1,\gamma}}{M_{s-m}},
\end{equation*}
so that
\begin{equation}
\label{rotxis_est_2}  
\sum_{s>0}\left\|\nabla^{\rm L} \times \xi^{(s)}\right\|_{0,\gamma}\, M_s^{-1}\, t^s
\leq  \frac{1}{M_1}\|\omega_0\|_{0,\gamma}t +\frac{3}{2}C_{a}M_0 \zeta^2.
\end{equation}
Multiplying (\ref{xisest_1}) by $M_s^{-1}\,t^s$,  summing the result over the index $s$,
and using  (\ref{divxis_est_2})--(\ref{rotxis_est_2}), we then obtain
\begin{equation}
\label{FEBB}
\zeta(t) \leq C_{\rm DN} \left(
\frac{15}{2}C_aM_0\zeta^2(t) +6C_a^2M_0^2\zeta^3(t) +  \frac{1}{M_1}\|\omega_0\|_{0,\gamma}t+
\sum_{s>0}\| \xi^{(s)}\cdot \nu\|_{1,\gamma}\, M_s^{-1}\, t^s
\right).
\end{equation}  
It remains to control the boundary term in the left-hand side of (\ref{FEBB}). Since we assumed
$\mathcal{S}\in  \mathcal{C}\{\mathcal{M}\}(\partial\Omega;\R)$, there exist two constants
$C_{\mathcal{S}}$ and $R_{\mathcal{S}}$ such that
\begin{equation}
\label{est_Surf}
\left\|\partial^\beta \mathcal{S} \right\|_{L^\infty} \leq C_{\mathcal{S}} R_{\mathcal{S}}^{-|\beta|}|\beta| !\, M_{|\beta|}.
\end{equation}  
Let us define
\begin{equation*}
\label{nota_supp_1}
C_{\nabla \mathcal{S}} :=  \left\| |\nabla \mathcal{S}(q)|_2^{-1} \right\|_{1,\gamma}. 
\end{equation*}
Then, using (\ref{NBC_1})--(\ref{def_nu_by_S}),  (\ref{algebra-Holder-m})  and (\ref{est_Surf}),
we obtain
\begin{eqnarray}
\label{NBT_1}
\frac{\| \xi^{(s)}\cdot \nu\|_{1,\gamma}}{M_s} &\leq &
C_aC_{\nabla \mathcal{S}} \| \xi^{(s)}\cdot \nabla\mathcal{S} \|_{1,\gamma} M_s^{-1}\nonumber \\
&\leq&C_{\nabla \mathcal{S}} \sum_{1\leq |\beta| \leq s}
\left(\frac{C_a}{R_{\mathcal{S}}}\right)^{|\beta|+1}\frac{C_{\mathcal{S}}}{R_{\mathcal{S}}}
(|\beta|+2)!\,\frac{M_{|\beta|+2}}{M_s}
\sum_{i=1}^s \sum_{P_i(s,\beta)} \prod_{j=1}^i 
\frac{\|\xi^{(\ell_j)}\|_{1,\gamma}^{|k_j|}}{k_{j}!} \nonumber\\
&\leq&C_{\nabla \mathcal{S}} \frac{C_{\mathcal{S}}}{R_{\mathcal{S}}}\sum_{1\leq |\beta| \leq s}
\left(\frac{C_a}{R_{\mathcal{S}}}\right)^{|\beta|+1} (|\beta|+2)(|\beta|+1)
|\beta|!\,\frac{M_{|\beta|+2}}{M_{|\beta|+1}}\frac{M_{|\beta|+1}}{M_{|\beta|}}\nonumber\\
&& \sum_{i=1}^s \sum_{P_i(s,\beta)}
\frac{M_{|\beta|}M_{\ell_1}^{|k_1|}\ldots M_{\ell_i}^{|k_i|}}{M_s}
\prod_{j=1}^i\left( \frac{\|\xi^{(\ell_j)}\|_{1,\gamma}}{ M_{\ell_j}}\right)^{|k_j|}\frac{1}{k_{j}!}.
\end{eqnarray}  
It is now convenient to use the notation
\[
\alpha_1:=\ell_1,\ \ldots,\ \,\, \alpha_{|k_1|}:=\ell_1,\ \,\, \alpha_{|k_1|+1}:=\ell_2,\ \ldots,\ \,\, 
\alpha_{|k_1|+|k_2|}:=\ell_2,\ \ldots,\ \,\, \alpha_{|k_1|+\ldots+|k_i|}:=\ell_i,
\]
in terms of which we have
\begin{eqnarray}
 M_{\ell_1}^{|k_1|}\,\ldots\, M_{\ell_i}^{|k_i|} &=& 
M_{\alpha_1}\,\ldots \,M_{\alpha_{|k_1|}} M_{\alpha_{|k_1|+1}}\,\ldots\, M_{\alpha_{|k_1|+|k_2|}}\,\ldots\, M_{\alpha_{|k_1|+\ldots+|k_i|}} \nonumber\\
&=& M_{\alpha_1}\, \ldots\, M_{\alpha_{|\beta|}}.
\label{index_recast}
\end{eqnarray}
Using (\ref{index_recast}) and (\ref{GenFunc}), equation \eqref{NBT_1} becomes
\begin{multline}
\frac{\| \xi^{(s)}\cdot \nu\|_{1,\gamma}}{M_s} \leq
C_a \frac{C_{\mathcal{S}}}{R_{\mathcal{S}}}\frac{C_{\nabla \mathcal{S}}}{R_{\mathcal{S}}}\sum_{1\leq |\beta| \leq s}
(|\beta|+2)(|\beta|+1)\left(\frac{C_a}{R_{\mathcal{S}}}\right)^{|\beta|} 
|\beta|!\,\frac{M_{|\beta|+2}}{M_{|\beta|+1}}\frac{M_{|\beta|+1}}{M_{|\beta|}}\\
\sum_{i=1}^s \sum_{P_i(s,\beta)}
\frac{M_{|\beta|}M_{\alpha_1} \ldots M_{\alpha_{|\beta|}}}{M_s}
\prod_{j=1}^i\frac{\left(\partial^{\ell_j}\zeta(0)\right)^{|k_j|}}{k_{j}! (\ell_j!)^{|k_j|}}.
\label{NBT_2}
\end{multline}
On the one hand, the differentiability property \eqref{Diff-P} leads to
\begin{equation}
\label{BoundDiff}
\frac{M_{|\beta|+2}}{M_{|\beta|+1}}\frac{M_{|\beta|+1}}{M_{|\beta|}} \leq C_{\rm d}^{2|\beta|+1}.
\end{equation}
On the other hand, the FdB-stability property \eqref{log-SL-P} leads to
\begin{equation}
\label{BoundFdB}
\frac{M_{|\beta|}M_{\alpha_1} \ldots M_{\alpha_{|\beta|}}}{M_s} \leq \frac{M_{\alpha_1 +\ldots+ \alpha_{|\beta|}}}{M_s} \leq  1.
\end{equation}
From \eqref{BoundDiff}--\eqref{BoundFdB}, estimate \eqref{NBT_2} becomes
\begin{eqnarray}
\frac{\| \xi^{(s)}\cdot \nu\|_{1,\gamma}}{M_s} &\leq&
C_a C_{\rm d} \frac{C_{\mathcal{S}}}{R_{\mathcal{S}}}\frac{C_{\nabla \mathcal{S}}}{R_{\mathcal{S}}}
\frac{(s+2)(s+1)}{s!} \nonumber\\
&& s!\sum_{1\leq |\beta| \leq s}
\left(\frac{C_a C_{\rm d}}{R_{\mathcal{S}}}\right)^{|\beta|} 
|\beta|! \sum_{i=1}^s \sum_{P_i(s,\beta)}
\prod_{j=1}^i\frac{\left(\partial^{\ell_j}\zeta(0)\right)^{|k_j|}}{k_{j}! (\ell_j!)^{|k_j|}} \nonumber\\
&\leq& C_a C_{\rm d}^2 \frac{C_{\mathcal{S}}}{R_{\mathcal{S}}}\frac{C_{\nabla \mathcal{S}}}{R_{\mathcal{S}}}
\frac{(s+2)(s+1)}{s!} \nonumber\\
&& s!\sum_{1\leq |\beta| \leq s} (\partial^\beta G)(0,0,0)
 \sum_{i=1}^s \sum_{P_i(s,\beta)}
\prod_{j=1}^i\frac{\left(\partial^{\ell_j}\zeta(0)\right)^{|k_j|}}{k_{j}! (\ell_j!)^{|k_j|}}, \label{almost_done}
\end{eqnarray}
where the map $G:\R^3\rightarrow \R$ is defined by 
\[
G(x_1,x_2,x_3) = \prod_{i=1}^3\left({1-\frac{C_a C_{\rm d}^2}{R_{\mathcal{S}}}}x_i \right)^{-1/3}, 
\]
and satisfies
\[
(\partial^\beta G)(0,0,0)= \partial_z^{|\beta|} G(z,z,z)_{|_{z=0}} = 
|\beta|! \left(\frac{C_a C_{\rm d}^2}{R_{\mathcal{S}}}\right)^{|\beta|}. 
\]
Setting
\[
g^{(s)}(0) := s!\sum_{1\leq |\beta| \leq s} (\partial^\beta G)(0,0,0)
 \sum_{i=1}^s \sum_{P_i(s,\beta)}
\prod_{j=1}^i\frac{\left(\partial^{\ell_j}\zeta(0)\right)^{|k_j|}}{k_{j}! (\ell_j!)^{|k_j|}}, 
\] 
\begin{equation}
\label{some_cst_def}
C_{{\rm d}a\mathcal{S}} := C_a C_{\rm d}^2 \frac{C_{\mathcal{S}}}{R_{\mathcal{S}}}\frac{C_{\nabla \mathcal{S}}}{R_{\mathcal{S}}},
\quad \quad 
C_{\mathcal{S}a{\rm d}}^{-1} := \frac{C_a C_{\rm d}^2}{R_{\mathcal{S}}},
\end{equation}
and using the Fa\`a di Bruno formula, we can rewrite \eqref{almost_done} as 
\begin{eqnarray}
\sum_{s> 0}{\| \xi^{(s)}\cdot \nu\|_{1,\gamma}}M_s^{-1}t^s &\leq&
C_{{\rm d}a\mathcal{S}}
\sum_{s>0} (s+2)(s+1) \frac{g^{(s)}(0)}{s!} t^s \nonumber \\
&\leq& C_{{\rm d}a\mathcal{S}}
\frac{d^2}{dt^2}\left(t^2g(t)\right) \nonumber \\
&\leq& C_{{\rm d}a\mathcal{S}}
\frac{d^2}{dt^2}\left(t^2G\left(\zeta(t),\zeta(t),\zeta(t)\right)\right) \nonumber \\
&\leq& C_{{\rm d}a\mathcal{S}}
\frac{d^2}{dt^2}\left(\frac{t^2}{1-C_{\mathcal{S}a{\rm d}}^{-1} \,\zeta(t)}\right). \label{final_bobo}
\end{eqnarray}
Combining \eqref{FEBB} and \eqref{final_bobo}, we obtain the final a priori estimate
\begin{equation*}
\label{final_ae_1}
\zeta(t) \leq C_{\rm DN} \left(
\frac{15}{2}C_aM_0\zeta^2(t) +6C_a^2M_0^2\zeta^3(t) +  \frac{1}{M_1}\|\omega_0\|_{0,\gamma}t+
 C_{{\rm d}a\mathcal{S}}
\frac{d^2}{dt^2}\left(\frac{t^2}{1-C_{\mathcal{S}a{\rm d}}^{-1} \,\zeta(t)}\right)
\right),
\end{equation*}
which we rewrite as
\begin{equation}
\label{final_ae_2}
\frac{d^2}{dt^2}\left(\frac{t^2}{C_{\mathcal{S}a{\rm d}}^{-1}
  \,\zeta(t)-1}\right) \leq h(t),
\end{equation}
where
\begin{eqnarray*}
  h(t)&:=&  C_{{\rm d}a\mathcal{S}}^{-1} 
\left( 6C_a^2M_0^2\zeta^3(t) + \frac{15}{2}C_aM_0\zeta^2(t) -C_{\rm DN}^{-1} \zeta(t) +  \frac{1}{M_1}\|\omega_0\|_{0,\gamma}t
 \right) \\[0.5ex]
&=&  \frac{d^2}{dt^2}\left( \int_0^t (t-\tau) h(\tau)d\tau\right).
\end{eqnarray*}
Twice integrating the differential inequality 
\eqref{final_ae_2}, and assuming that \hbox{$1-C_{\mathcal{S}a{\rm d}}^{-1}
  \,\zeta(t)\geq 0$,} for the time being, we
obtain
\begin{equation*}
\label{final_ae_3}
-1 \leq \left(1-C_{\mathcal{S}a{\rm d}}^{-1} \,\zeta(t)\right) \frac{1}{t^2} \int_0^th(\tau)d\tau.
\end{equation*}
A sufficient condition for this 
differential inequality to hold is to have simultaneously
\begin{equation}
\label{SufCond}
1-C_{\mathcal{S}a{\rm d}}^{-1} \,\zeta(t)\geq 0, \quad {\rm and} \quad 
Q(\zeta):=6C_a^2M_0^2\zeta^3(t) + \frac{15}{2}C_aM_0\zeta^2(t) -C_{\rm DN}^{-1} \zeta(t) + \Gamma(t) \geq 0,
\end{equation}
where we have set $\Gamma(t):= (\|\omega_0\|_{0,\gamma}/M_1)\,t $.
The discriminant of the cubic polynomial $Q(\zeta)$,
\[
\Delta := \frac{1}{54C_a^6M_0^6}\left(\frac{25}{8} +\frac{1}{C_{\rm DN}} \right)^3
-\frac{9}{C_a^4M_0^4}\left(\frac{5}{72}\frac{1}{C_aM_0}\left(\frac{25}{12} +\frac{1}{C_{\rm DN}}\right) + \frac{\Gamma}{6}\right)^2,
\]
is positive at small times, when $\Gamma$ is sufficiently small, whereby $Q(\zeta)$ has three
real roots $\zeta_i$. The polynomial has two local extrema at points of different signs, hence
it has roots of both signs. Since by Vi\`ete's theorem the product of the roots is negative,
two roots are positive and one is negative. The second inequality of \eqref{SufCond} implies
\begin{equation}
\label{bound_zeta}
  \zeta \leq \zeta_2(\Gamma),
\end{equation}
where $\zeta_2$ is the intermediate root, i.e. the smaller of the two positive roots.

We determine now the largest time $t_c$, for which bound \eqref{bound_zeta} holds.
Differentiating the equation for the roots of polynomial \eqref{SufCond}
with respect to $\Gamma$,
we find
\[
\frac{\partial \zeta_i}{\partial \Gamma}=-\left({\frac{\partial Q}{\partial \zeta}} \Bigg |_{\zeta=\zeta_i} \right)^{-1}.
\]
Consequently, on increasing $\Gamma$ from zero, the root $\zeta_2$ monotonically increases from zero and $\zeta_3$
monotonically decreases till the two roots collide (i.e. $\zeta_2=\zeta_3$) at the critical value
$\Gamma=\Gamma_c$, determined by the equation $\Delta=0$.
This corresponds to the critical time $t_c=\Gamma_cM_1/\|\omega_0\|_{0,\gamma}$.
Setting $t_{\mathcal{S}a{\rm d}}=Q(C_{\mathcal{S}a{\rm d}})M_1/\|\omega_0\|_{0,\gamma}$, an upper bound of the radius of convergence $T$
of the generating function $\zeta$ is given by the smaller of $t_{\mathcal{S}a{\rm d}}$ and $t_c$. Therefore, we obtain
\[
\zeta(t) \leq \min\{C_{\mathcal{S}a{\rm d}},\, \zeta_2(\Gamma_c) \}, \quad t\in]-T,T[, \quad T=\min\{t_{\mathcal{S}a{\rm d}}, t_c\},
\]
and the sufficient condition \eqref{SufCond} is indeed satisfied. We observe, from
the definition of $t_{\mathcal{S}a{\rm d}}$ and $t_c$ that the maximum time for which the
Lagrangian map is time-ultradifferentiable is controlled by the H\"older
continuity of the initial vorticity $\|\omega_0\|_{0,\gamma}$ and the 
``radius of ultradifferentiability'' $R_{\mathcal{S}}$ of the boundary.

Finally, we note that our solution was constructed from a Taylor expansion around $t=0$. However, we 
can restart the process from any time $t_\ast>0$, at which the vorticity has $\mathcal{C}^{0,\eta}$-regularity
for some $\eta >0$ (which may be less than $\gamma$). Then,  ultradifferentiability in time 
persists on the interval $[t_\ast,t_\ast + T_\ast[$ for some positive $T_\ast$. This 
process can be continued until the vorticity ceases to have any $\mathcal{C}^{0,\eta}$-regularity.
Clearly the Lagrangian map remains time-ultradifferentiable as long as the vorticity 
is H\"older continuous. This completes our proof.

\section{Concluding remarks}
\label{s:concl}

For 3D incompressible Euler flow with an initial vorticity that is
H\"older continuous, we have extended the results of Refs.~\cite{FZ14, ZF14} on time-analyticity
of the Lagrangian map to the case of wall-bounded flow with a a fixed
impermeable  boundary  that is analytic in its shape. We have also obtained 
similar regularity results when the boundary is in a suitable class of
ultradifferentiability. 

An important feature that singles out our proof is that it is constructive, in
the sense that it leads to explicit recursion relations for the time-Taylor
coefficients of the Lagrangian map. This feature was crucial in allowing the
development of an efficient Cauchy-Lagrangian numerical scheme for the case
without boundary \cite{PZF16}. Such a scheme can in principle be carrried over
to the case with a boundary; this will, of course, require the numerical
solution of the elliptic (Poisson) equations with both Dirichlet and Neumann
boundary conditions involved in the Helmholtz--Hodge decompositions. For the
case of a boundary that is not analytic, but in a suitable 
ultradifferential class, the time-Taylor series will generally be
divergent for any $t>0$, so that some resummation method (e.g., 
by a Borel transformation  \cite{PF07}) must also be used for the numerical
implementation.

Once a variant of the Cauchy-Lagrangian scheme, adapted to
solid boundaries has been developped, it will be of interest to revisit
the issue of finite time blow-up via state-of-the-art simulations.

So far, we have only considered the case of flow within a bounded simply
connected domain. The extension of our results to outer flow (e.g. 
flow in an infinite domain around a bounded obstacle) and to multiply
connected domains will be considered later. 

It is also of interest to handle the extension of our results to other flow 
equations. One obvious case is the Euler--Poisson equation for which a
Lagrangian time-analyticity result was obtained \cite{ZF14, RVF15}. This was done
in the absence of boundary and in a cosmological context of attractive
gravitational forces. Boundaries are usually not considered in cosmology.
However, they play an important part in plasma physics because, for example,
fusion facilities are always enclosed. Of particular interest for plasmas
is the  Euler--Poisson equation with a repulsive electrostatic force,
which will be considered in subsequent work. 

\section*{Acknowledgements}
We are grateful to Claude Bardos and L\'aszl\'o  Sz\'ekelyhidi for many fruitful discussions.
This work was supported by the VLASIX and EUROFUSION projects respectively
under the grants No ANR-13-MONU-0003-01 and EURATOM-WP15-ENR-01/IPP-01. We
have also benefitted from our participation in the Wolfgang Pauli Institute
workshop ``Euler and Navier-Stokes Equations and Connected Topics,'' Vienna, December
14-18, 2015.

\end{document}